\documentclass[12pt,a4paper]{article}
\usepackage{amsfonts}

\usepackage{graphicx}
\usepackage{amsmath}

\newtheorem{theorem}{Theorem}

\newtheorem{lemma}[theorem]{Lemma}

\newtheorem{proposition}[theorem]{Proposition}
\newtheorem{remark}[theorem]{Remark}

\newenvironment{proof}[1][Proof]{\noindent\textbf{#1.} }{\hfill $\Box $}

\begin{document}

\title{Weak analytic hyperbolicity of generic hypersurfaces of high degree in $%
\mathbb{P}^{4} $}
\author{Erwan Rousseau \\
{\small {D\'epartement de Math\'ematiques,} }\\
{\small {Universit\'e du Qu\'ebec \`a Montr\'eal,} }\\
{\small {(e-mail: eroussea@math.uqam.ca)}}}
\date{}
\maketitle

\begin{abstract}
In this article we prove that every entire curve in a generic hypersurface
of degree $d\geq 593$ in $\mathbb{P}_{\mathbb{C}}^{4}$ is algebraically
degenerated i.e there exists a proper subvariety which contains the entire
curve. \newline
\begin{center}
\textbf{R\'esum\'e}
\end{center}
Dans cet article nous d\'emontrons que toute courbe enti\`ere dans
une
hypersurface g\'en\'erique de degr\'e $d\geq 593$ dans $\mathbb{P}_{\mathbb{C%
}}^{4}$ est alg\'ebriquement d\'eg\'en\'er\'ee i.e il existe une
sous-vari\'et\'e propre qui contient la courbe enti\`ere.
\end{abstract}

\section{Introduction}

In 1970, S. Kobayashi conjectured in \cite{Ko70} that a generic hypersurface
$X$ in $\mathbb{P}_{\mathbb{C}}^{n}$ is hyperbolic provided that $d=\deg
(X)\geq 2n-1,$ for $n\geq 3.$ For $n=3,$ it was obtained by Demailly and El
Goul in \cite{DEG00} that $d\geq 21$ implies the hyperbolicity of very
generic hypersurfaces $X$ in $\mathbb{P}_{\mathbb{C}}^{3}.$

In this paper, we would like to prove the following:

\begin{theorem}
Let $X\subset \mathbb{P}_{\mathbb{C}}^{4}$ be a generic hypersurface such
that $d=deg(X)\geq 593.$ Then every entire curve $f:\mathbb{C}\rightarrow X$
is algebraically degenerated, i.e there exists a proper subvariety $Y\subset
X$ such that $f(\mathbb{C)}\subset Y.$
\end{theorem}

This result is a weaker version of the conjecture in dimension 3, because to
obtain the full conjecture one needs to prove that the entire curve is
constant.

The proof of the theorem is based on two techniques. Consider $\mathcal{X}%
\subset \mathbb{P}^{4}\times \mathbb{P}^{N_{d}}$ the universal hypersurface
of degree $d$ in $\mathbb{P}^{4}.$

In the first section we construct meromorphic vector fields on the space $%
J_{3}^{v}(\mathcal{X})$ of vertical 3-jets of $\mathcal{X}.$ This technique,
initiated by Clemens \cite{Cle}, Ein \cite{Ein}, Voisin \cite{Voi}, was
generalized by Y.-T. Siu \cite{SY04} and detailed by M. Paun in dimension 2
\cite{Paun05}. Here we generalize M. Paun's computations in dimension 3 and
obtain that a pole order equal to 12 is enough to obtain a ''large'' space
of global sections of the twisted tangent bundle.

In the second section, we summarize the main facts about the bundle of jet
differentials of order $k$ and degree $m,$ $E_{k,m}T_{X}^{\ast }.$ The idea,
in hyperbolicity questions, is that global sections of this bundle vanishing
on ample divisors provide algebraic differential equations for any entire
curve $f:\mathbb{C}\rightarrow X.$ Therefore, the main point is to produce
enough algebraically independent global holomorphic jet differentials. In
the case of surfaces in $\mathbb{P}^{3}$ of degree $d\geq 15,$ one can
produce global jet differentials of order 2 vanishing on ample divisors
using a Riemann-Roch computation and a Bogomolov vanishing theorem (see \cite
{De95}). Y.-T. Siu, in \cite{SY04}, described a way to produce global jet
differentials vanishing on ample divisors for hypersurfaces of sufficiently
large degree $d$ in $\mathbb{P}^{n},$ for any $n.$ One problem is that the
bound obtained for $d$ is quite high. If we are interested in the degree $d$
for smooth hypersurfaces of\textit{\ }$\mathbb{P}^{4}$, an interesting
result obtained in \cite{Rou2}, is the existence of global jet differentials
of order 3 vanishing on ample divisors for $d\geq 97.$

In the last section we complete the proof of the theorem using Siu's
approach \cite{SY04}, by taking the derivative of the jet differential in
the direction of the vector fields constructed in the first part, and Mihai
Paun's ''trick'' described in \cite{Paun05} to avoid the use of McQuillan
results (cf. \cite{mcq}).

\bigskip

\textbf{Acknowledgements. }We would like to thank Mihai Paun for his
lectures about Siu's ideas given at the summer school Pragmatic in Catania,
2004.

\section{Vector fields}

In this first section, we generalize to dimension 3 the approach of Mihai
Paun \cite{Paun05} which gives some precisions to Siu's ideas \cite{SY04} in
dimension 2. Consider $\mathcal{X}\subset \mathbb{P}^{4}\times \mathbb{P}%
^{N_{d}}$ the universal hypersurface given by the equation
\begin{equation*}
\underset{\left| \alpha \right| =d}{\sum }a_{\alpha }Z^{\alpha }=0,\text{
where }[a]\in \mathbb{P}^{N_{d}}\text{ and }[Z]\in \mathbb{P}^{4}.
\end{equation*}

We use the notations: for $\alpha =(\alpha _{0},...,\alpha _{4})\in \mathbb{N%
}^{5},$ $\left| \alpha \right| =\sum_{i}\alpha _{i}$ and if $%
Z=(Z_{0},Z_{1},...,Z_{4})$ are homogeneous coordinates on $\mathbb{P}^{4},$
then $Z^{\alpha }=\prod Z_{j}^{\alpha _{j}}.$ $\mathcal{X}$ is a smooth
hypersurface of degree $(d,1)$ in $\mathbb{P}^{4}\times \mathbb{P}^{N_{d}}.$
We denote by $J_{3}(\mathcal{X})$ the manifold of the 3-jets in $\mathcal{X}%
, $ and $J_{3}^{v}(\mathcal{X})$ the submanifold of $J_{3}(\mathcal{X})$
consisting of 3-jets in $\mathcal{X}$ tangent to the fibers of the
projection $\pi :\mathcal{X}\rightarrow \mathbb{P}^{N_{d}}.$

Let us consider the set $\Omega _{0}:=(Z_{0}\neq 0)\times (a_{0d000}\neq
0)\subset \mathbb{P}^{4}\times \mathbb{P}^{N_{d}}.$ We assume that global
coordinates are given on $\mathbb{C}^{4}$ and $\mathbb{C}^{N_{d}}.$ The
equation of $\mathcal{X}$ becomes
\begin{equation*}
\mathcal{X}_{0}:=(z_{1}^{d}+\underset{\left| \alpha \right| \leq d,\text{ }%
\alpha _{1}<d}{\sum }a_{\alpha }z^{\alpha }=0)
\end{equation*}

Then the equations of $J_{3}^{v}(\mathcal{X}_{0})$ in $\mathbb{C}^{4}\times
\mathbb{C}^{N_{d}}\times \mathbb{C}^{4}\times \mathbb{C}^{4}\times \mathbb{C}%
^{4}$ can be written:
\begin{equation}
\underset{\left| \alpha \right| \leq d,\text{ }a_{d000}=1}{\sum }a_{\alpha
}z^{\alpha }=0
\end{equation}
\begin{equation}
\underset{j=1}{\overset{4}{\sum }}\underset{\left| \alpha \right| \leq d,%
\text{ }a_{d000}=1}{\sum }a_{\alpha }\frac{\partial z^{\alpha }}{\partial
z_{j}}\xi _{j}^{(1)}=0
\end{equation}
\begin{equation}
\underset{j=1}{\overset{4}{\sum }}\underset{\left| \alpha \right| \leq d,%
\text{ }a_{d000}=1}{\sum }a_{\alpha }\frac{\partial z^{\alpha }}{\partial
z_{j}}\xi _{j}^{(2)}+\underset{j,k=1}{\overset{4}{\sum }}\underset{\left|
\alpha \right| \leq d,\text{ }a_{d000}=1}{\sum }a_{\alpha }\frac{\partial
^{2}z^{\alpha }}{\partial z_{j}\partial z_{k}}\xi _{j}^{(1)}\xi _{k}^{(1)}=0
\end{equation}
\begin{eqnarray}
\underset{j=1}{\overset{4}{\sum }}\underset{\left| \alpha \right| \leq d,%
\text{ }a_{d000}=1}{\sum }a_{\alpha }\frac{\partial z^{\alpha }}{\partial
z_{j}}\xi _{j}^{(3)}+\underset{j,k=1}{\overset{4}{3\sum }}\underset{\left|
\alpha \right| \leq d,\text{ }a_{d000}=1}{\sum }a_{\alpha }\frac{\partial
^{2}z^{\alpha }}{\partial z_{j}\partial z_{k}}\xi _{j}^{(2)}\xi _{k}^{(1)} &&
\notag \\
+\underset{j,k,l=1}{\overset{4}{\sum }}\underset{\left| \alpha \right| \leq
d,\text{ }a_{d000}=1}{\sum }a_{\alpha }\frac{\partial ^{3}z^{\alpha }}{%
\partial z_{j}\partial z_{k}\partial z_{l}}\xi _{j}^{(1)}\xi _{k}^{(1)}\xi
_{l}^{(1)} &=&0
\end{eqnarray}

Consider now a vector field
\begin{equation*}
V=\underset{\left| \alpha \right| \leq d,\text{ }\alpha _{1}<d}{\sum }%
v_{\alpha }\frac{\partial }{\partial a_{\alpha }}+\underset{j}{\sum }v_{j}%
\frac{\partial }{\partial z_{j}}+\underset{j,k}{\sum }w_{j,k}\frac{\partial
}{\partial \xi _{j}^{(k)}}
\end{equation*}
on the vector space $\mathbb{C}^{4}\times \mathbb{C}^{N_{d}}\times \mathbb{C}%
^{4}\times \mathbb{C}^{4}\times \mathbb{C}^{4}.$ The conditions to be
satisfied by $V$ to be tangent to $J_{3}^{v}(\mathcal{X}_{0})$ are
\begin{equation*}
\underset{\left| \alpha \right| \leq d,\text{ }\alpha _{1}<d}{\sum }%
v_{\alpha }z^{\alpha }+\underset{j=1}{\overset{4}{\sum }}\underset{\left|
\alpha \right| \leq d,\text{ }a_{d000}=1}{\sum }a_{\alpha }\frac{\partial
z^{\alpha }}{\partial z_{j}}v_{j}=0
\end{equation*}
\begin{equation*}
\underset{j=1}{\overset{4}{\sum }}\underset{\left| \alpha \right| \leq d,%
\text{ }\alpha _{1}<d}{\sum }v_{\alpha }\frac{\partial z^{\alpha }}{\partial
z_{j}}\xi _{j}^{(1)}+\underset{j,k=1}{\overset{4}{\sum }}\underset{\left|
\alpha \right| \leq d,\text{ }a_{d000}=1}{\sum }a_{\alpha }\frac{\partial
^{2}z^{\alpha }}{\partial z_{j}\partial z_{k}}v_{j}\xi _{k}^{(1)}+\underset{%
j=1}{\overset{4}{\sum }}\underset{\left| \alpha \right| \leq d,\text{ }%
a_{d000}=1}{\sum }a_{\alpha }\frac{\partial z^{\alpha }}{\partial z_{j}}%
w_{j}^{(1)}=0
\end{equation*}
\begin{eqnarray*}
\underset{\left| \alpha \right| \leq d,\text{ }\alpha _{1}<d}{\sum }(%
\underset{j=1}{\overset{4}{\sum }}\frac{\partial z^{\alpha }}{\partial z_{j}}%
\xi _{j}^{(2)}+\underset{j,k=1}{\overset{4}{\sum }}\frac{\partial
^{2}z^{\alpha }}{\partial z_{j}\partial z_{k}}\xi _{j}^{(1)}\xi
_{k}^{(1)})v_{\alpha } && \\
+\underset{j=1}{\overset{4}{\sum }}\underset{\left| \alpha \right| \leq d,%
\text{ }a_{d000}=1}{\sum }a_{\alpha }(\underset{k=1}{\overset{4}{\sum }}%
\frac{\partial ^{2}z^{\alpha }}{\partial z_{j}\partial z_{k}}\xi _{k}^{(2)}+%
\underset{k,l=1}{\overset{4}{\sum }}\frac{\partial ^{3}z^{\alpha }}{\partial
z_{j}\partial z_{k}\partial z_{l}}\xi _{k}^{(1)}\xi _{l}^{(1)})v_{j} && \\
+\underset{\left| \alpha \right| \leq d,\text{ }a_{d000}=1}{\sum }(\underset{%
j,k=1}{\overset{4}{\sum }}a_{\alpha }\frac{\partial ^{2}z^{\alpha }}{%
\partial z_{j}\partial z_{k}}(w_{j}^{(1)}\xi _{k}^{(1)}+w_{k}^{(1)}\xi
_{j}^{(1)})+\underset{j=1}{\overset{4}{\sum }}a_{\alpha }\frac{\partial
z^{\alpha }}{\partial z_{j}}w_{j}^{(2)}) &=&0
\end{eqnarray*}
\begin{eqnarray*}
\underset{\left| \alpha \right| \leq d,\text{ }\alpha _{1}<d}{\sum }(%
\underset{j=1}{\overset{4}{\sum }}\frac{\partial z^{\alpha }}{\partial z_{j}}%
\xi _{j}^{(3)}+3\underset{j,k=1}{\overset{4}{\sum }}\frac{\partial
^{2}z^{\alpha }}{\partial z_{j}\partial z_{k}}\xi _{j}^{(2)}\xi _{k}^{(1)}+%
\underset{j,k,l=1}{\overset{4}{\sum }}\frac{\partial ^{3}z^{\alpha }}{%
\partial z_{j}\partial z_{k}\partial z_{l}}\xi _{j}^{(1)}\xi _{k}^{(1)}\xi
_{l}^{(1)})v_{\alpha } && \\
+\underset{j=1}{\overset{4}{\sum }}\underset{\left| \alpha \right| \leq d,%
\text{ }a_{d000}=1}{\sum }a_{\alpha }(\underset{k=1}{\overset{4}{\sum }}%
\frac{\partial ^{2}z^{\alpha }}{\partial z_{j}\partial z_{k}}\xi _{k}^{(3)}+3%
\underset{k,l=1}{\overset{4}{\sum }}\frac{\partial ^{3}z^{\alpha }}{\partial
z_{j}\partial z_{k}\partial z_{l}}\xi _{k}^{(2)}\xi _{l}^{(1)}+\underset{%
k,l,m=1}{\overset{4}{\sum }}\frac{\partial ^{4}z^{\alpha }}{\partial
z_{j}\partial z_{k}\partial z_{l}\partial z_{m}}\xi _{k}^{(1)}\xi
_{l}^{(1)}\xi _{m}^{(1)})v_{j} && \\
+\underset{\left| \alpha \right| \leq d,\text{ }a_{d000}=1}{\sum }(\underset{%
j,k,l=1}{\overset{4}{\sum }}a_{\alpha }\frac{\partial ^{3}z^{\alpha }}{%
\partial z_{j}\partial z_{k}\partial z_{l}}(w_{j}^{(1)}\xi _{k}^{(1)}\xi
_{l}^{(1)}+\xi _{j}^{(1)}w_{k}^{(1)}\xi _{l}^{(1)}+\xi _{j}^{(1)}\xi
_{k}^{(1)}w_{l}^{(1)}) && \\
+3\underset{j,k=1}{\overset{4}{\sum }}a_{\alpha }\frac{\partial
^{2}z^{\alpha }}{\partial z_{j}\partial z_{k}}(w_{j}^{(2)}\xi _{k}^{(1)}+\xi
_{j}^{(2)}w_{k}^{(1)})+\underset{j=1}{\overset{4}{\sum }}a_{\alpha }\frac{%
\partial z^{\alpha }}{\partial z_{j}}w_{j}^{(3)})=0 &&
\end{eqnarray*}

Now we can introduce the first package of vector fields tangent to $%
J_{3}^{v}(\mathcal{X}_{0}).$ We denote by $\delta _{j}\in \mathbb{N}^{4}$
the multi-index whose j-component is equal to 1 and the other are zero.

For $\alpha _{1}\geq 4:$%
\begin{equation*}
V_{\alpha }^{4000}:=\frac{\partial }{\partial a_{\alpha }}-4z_{1}\frac{%
\partial }{\partial a_{\alpha -\delta _{1}}}+6z_{1}^{2}\frac{\partial }{%
\partial a_{\alpha -2\delta _{1}}}-4z_{1}^{3}\frac{\partial }{\partial
a_{\alpha -3\delta _{1}}}+z_{1}^{4}\frac{\partial }{\partial a_{\alpha
-4\delta _{1}}}.
\end{equation*}

For $\alpha _{1}\geq 3,\alpha _{2}\geq 1:$%
\begin{eqnarray*}
V_{\alpha }^{3100} &:&=\frac{\partial }{\partial a_{\alpha }}-3z_{1}\frac{%
\partial }{\partial a_{\alpha -\delta _{1}}}-z_{2}\frac{\partial }{\partial
a_{\alpha -\delta _{2}}}+3z_{1}z_{2}\frac{\partial }{\partial a_{\alpha
-\delta _{1}-\delta _{2}}} \\
&&+3z_{1}^{2}\frac{\partial }{\partial a_{\alpha -2\delta _{1}}}%
-3z_{1}^{2}z_{2}\frac{\partial }{\partial a_{\alpha -2\delta _{1}-\delta
_{2}}}-z_{1}^{3}\frac{\partial }{\partial a_{\alpha -3\delta _{1}}}%
+z_{1}^{3}z_{2}\frac{\partial }{\partial a_{\alpha -3\delta _{1}-\delta _{2}}%
}.
\end{eqnarray*}

For $\alpha _{1}\geq 2,\alpha _{2}\geq 2:$%
\begin{equation*}
V_{\alpha }^{2200}:=\frac{\partial }{\partial a_{\alpha }}-z_{2}\frac{%
\partial }{\partial a_{\alpha -\delta _{2}}}-z_{1}\frac{\partial }{\partial
a_{\alpha -\delta _{1}}}+z_{1}z_{2}^{2}\frac{\partial }{\partial a_{\alpha
-\delta _{1}-2\delta _{2}}}+z_{1}^{2}z_{2}\frac{\partial }{\partial
a_{\alpha -2\delta _{1}-\delta _{2}}}-z_{1}^{2}z_{2}^{2}\frac{\partial }{%
\partial a_{\alpha -2\delta _{1}-2\delta _{2}}}.
\end{equation*}

For $\alpha _{1}\geq 2,\alpha _{2}\geq 1,\alpha _{3}\geq 1:$%
\begin{eqnarray*}
V_{\alpha }^{2110} &:&=\frac{\partial }{\partial a_{\alpha }}-z_{3}\frac{%
\partial }{\partial a_{\alpha -\delta _{3}}}-z_{2}\frac{\partial }{\partial
a_{\alpha -\delta _{2}}}-2z_{1}\frac{\partial }{\partial a_{\alpha -\delta
_{1}}}+z_{2}z_{3}\frac{\partial }{\partial a_{\alpha -\delta _{2}-\delta
_{3}}} \\
&&+2z_{1}z_{3}\frac{\partial }{\partial a_{\alpha -\delta _{1}-\delta _{3}}}%
+2z_{1}z_{2}\frac{\partial }{\partial a_{\alpha -\delta _{1}-\delta _{2}}}%
+z_{1}^{2}\frac{\partial }{\partial a_{\alpha -2\delta _{1}}} \\
&&-2z_{1}z_{2}z_{3}\frac{\partial }{\partial a_{\alpha -\delta _{1}-\delta
_{2}-\delta _{3}}}-z_{1}^{2}z_{3}\frac{\partial }{\partial a_{\alpha
-2\delta _{1}-\delta _{3}}}-z_{1}^{2}z_{2}\frac{\partial }{\partial
a_{\alpha -2\delta _{1}-\delta _{2}}}+z_{1}^{2}z_{2}z_{3}\frac{\partial }{%
\partial a_{\alpha -2\delta _{1}-\delta _{2}-\delta _{3}}}.
\end{eqnarray*}

For $\alpha _{1}\geq 1,\alpha _{2}\geq 1,\alpha _{3}\geq 1,\alpha _{4}\geq
1: $%
\begin{eqnarray*}
V_{\alpha }^{1111} &:&=\frac{\partial }{\partial a_{\alpha }}-z_{1}\frac{%
\partial }{\partial a_{\alpha -\delta _{1}}}-z_{2}\frac{\partial }{\partial
a_{\alpha -\delta _{2}}}-z_{3}\frac{\partial }{\partial a_{\alpha -\delta
_{3}}} \\
&&-z_{4}\frac{\partial }{\partial a_{\alpha -\delta _{4}}}+z_{2}z_{3}z_{4}%
\frac{\partial }{\partial a_{\alpha -\delta _{2}-\delta _{3}-\delta _{4}}}%
+z_{1}z_{3}z_{4}\frac{\partial }{\partial a_{\alpha -\delta _{1}-\delta
_{3}-\delta _{4}}} \\
&&+z_{1}z_{2}z_{3}\frac{\partial }{\partial a_{\alpha -\delta _{1}-\delta
_{2}-\delta _{3}}}+z_{1}z_{2}z_{4}\frac{\partial }{\partial a_{\alpha
-\delta _{1}-\delta _{2}-\delta _{4}}}-z_{1}z_{2}z_{3}z_{4}\frac{\partial }{%
\partial a_{\alpha -\delta _{1}-\delta _{2}-\delta _{3}-\delta _{4}}}.
\end{eqnarray*}

Similar vector fields are constructed by permuting the z-variables, and
changing the index $\alpha $ as indicated by the permutation. The pole order
of the previous vector fields is equal to 4.

\bigskip

\begin{lemma}
For any $(v_{i})_{1\leq i\leq 4}\in \mathbb{C}^{4},$ there exist $v_{\alpha
}(a),$ with degree at most 1 in the variables $(a_{\gamma }),$ such that $V:=%
\underset{\alpha }{\sum }v_{\alpha }(a)\frac{\partial }{\partial a_{\alpha }}%
+\underset{j}{\sum }v_{j}\frac{\partial }{\partial z_{j}}$ is tangent to $%
J_{3}^{v}(\mathcal{X}_{0})$ at each point.
\end{lemma}

\begin{proof}
We impose the additional conditions of vanishing for the coefficients of $%
\xi _{j}^{(1)}$ in the second equation (respectively of $\xi _{j}^{(1)}\xi
_{k}^{(1)}$ in the third equation and $\xi _{j}^{(1)}\xi _{k}^{(1)}\xi
_{l}^{(1)}$ in the fourth equation) for any $1\leq j\leq k\leq l\leq 4$.
Then the coefficients of $\xi _{j}^{(2)}$ (respectively $\xi _{j}^{(2)}\xi
_{k}^{(1)}$ and $\xi _{j}^{(3)})$ are automatically zero in the third
(respectively fourth) equation. The resulting 35 equations are
\begin{equation*}
\underset{\left| \alpha \right| \leq d,\text{ }\alpha _{1}<d}{\sum }%
v_{\alpha }z^{\alpha }+\underset{j=1}{\overset{4}{\sum }}\underset{\left|
\alpha \right| \leq d,\text{ }a_{d000}=1}{\sum }a_{\alpha }\frac{\partial
z^{\alpha }}{\partial z_{j}}v_{j}=0
\end{equation*}

\begin{equation*}
\underset{\left| \alpha \right| \leq d,\text{ }\alpha _{1}<d}{\sum }%
v_{\alpha }\frac{\partial z^{\alpha }}{\partial z_{j}}+\underset{k=1}{%
\overset{4}{\sum }}\underset{\left| \alpha \right| \leq d,\text{ }a_{d000}=1%
}{\sum }a_{\alpha }\frac{\partial ^{2}z^{\alpha }}{\partial z_{j}\partial
z_{k}}v_{k}=0
\end{equation*}

\begin{equation*}
\underset{\alpha }{\sum }\frac{\partial ^{2}z^{\alpha }}{\partial
z_{j}\partial z_{k}}v_{\alpha }+\underset{l=1}{\overset{4}{\sum }}\underset{%
\left| \alpha \right| \leq d,\text{ }a_{d000}=1}{\sum }a_{\alpha }\frac{%
\partial ^{3}z^{\alpha }}{\partial z_{j}\partial z_{k}\partial z_{l}}v_{l}=0
\end{equation*}

\begin{equation*}
\underset{\alpha }{\sum }\frac{\partial ^{3}z^{\alpha }}{\partial
z_{j}\partial z_{k}\partial z_{l}}v_{\alpha }+\underset{m=1}{\overset{4}{%
\sum }}\underset{\left| \alpha \right| \leq d,\text{ }a_{d000}=1}{\sum }%
a_{\alpha }\frac{\partial ^{4}z^{\alpha }}{\partial z_{j}\partial
z_{k}\partial z_{l}\partial z_{m}}v_{m}=0
\end{equation*}

Now we can observe that if the $v_{\alpha }(a)$ satisfy the first equation,
they automatically satisfy the other ones because the $v_{\alpha }$ are
constants with respect to $z$. Therefore it is sufficient to find $%
(v_{\alpha })$ satisfying the first equation. We identify the coefficients
of $z^{\rho }=z_{1}^{\rho _{1}}$ $z_{2}^{\rho _{2}}$ $z_{3}^{\rho _{3}}$ $%
z_{4}^{\rho _{4}}:$%
\begin{equation*}
v_{\rho }+\underset{j=1}{\overset{4}{\sum }}a_{\rho +\delta _{j}}v_{j}(\rho
_{j}+1)=0
\end{equation*}
\end{proof}

Another family of vector fields can be obtained thanks to the generalization
to dimension 3 of a lemma (cf. \cite{Paun05}) given by Mihai Paun in
dimension 2. Consider a $4\times 4$-matrix $A=(A_{j}^{k})\in \mathcal{M}_{4}(%
\mathbb{C})$ and let $\widetilde{V}:=\underset{j,k}{\sum }w_{j}^{(k)}\frac{%
\partial }{\partial \xi _{j}^{(k)}},$ where $w^{(k)}:=A\xi ^{(k)},$ for $%
k=1,2,3.$

\begin{lemma}
There exist polynomials $v_{\alpha }(z,a):=\underset{\left| \beta \right|
\leq 3}{\sum }v_{\beta }^{\alpha }(a)z^{\beta }$ where each coefficient $%
v_{\beta }^{\alpha }$ has degree at most 1 in the variables $(a_{\gamma })$
such that
\begin{equation*}
V:=\underset{\alpha }{\sum }v_{\alpha }(z,a)\frac{\partial }{\partial
a_{\alpha }}+\widetilde{V}
\end{equation*}
is tangent to $J_{3}^{v}(\mathcal{X}_{0})$ at each point.
\end{lemma}

\begin{proof}
We impose the additional conditions of vanishing for the coefficients of $%
\xi _{j}^{(1)}$ in the second equation (respectively of $\xi _{j}^{(1)}\xi
_{k}^{(1)}$ in the third equation and $\xi _{j}^{(1)}\xi _{k}^{(1)}\xi
_{l}^{(1)}$ in the fourth equation) for any $1\leq j\leq k\leq l\leq 4$.
Then the coefficients of $\xi _{j}^{(2)}$ (respectively $\xi _{j}^{(2)}\xi
_{k}^{(1)}$ and $\xi _{j}^{(3)})$ are automatically zero in the third
(respectively fourth) equation. The resulting 35 equations are
\begin{equation*}
\underset{\left| \alpha \right| \leq d,\text{ }\alpha _{1}<d}{\sum }%
v_{\alpha }z^{\alpha }=0\text{ \ }(5)
\end{equation*}

\begin{equation*}
\underset{\left| \alpha \right| \leq d,\text{ }\alpha _{1}<d}{\sum }%
v_{\alpha }\frac{\partial z^{\alpha }}{\partial z_{j}}+\underset{k=1}{%
\overset{4}{\sum }}\underset{\left| \alpha \right| \leq d,\text{ }a_{d000}=1%
}{\sum }a_{\alpha }\frac{\partial z^{\alpha }}{\partial z_{k}}A_{k}^{j}=0%
\text{ \ \ }(6_{j})
\end{equation*}

\begin{equation*}
\underset{\alpha }{\sum }\frac{\partial ^{2}z^{\alpha }}{\partial
z_{j}\partial z_{k}}v_{\alpha }+\underset{\alpha ,p}{\sum }a_{\alpha }\frac{%
\partial ^{2}z^{\alpha }}{\partial z_{j}\partial z_{p}}A_{p}^{k}+\underset{%
\alpha ,p}{\sum }a_{\alpha }\frac{\partial ^{2}z^{\alpha }}{\partial
z_{k}\partial z_{p}}A_{p}^{j}=0\text{ \ }(7_{jk})
\end{equation*}

\begin{equation*}
\underset{\alpha }{\sum }\frac{\partial ^{3}z^{\alpha }}{\partial
z_{j}\partial z_{k}\partial z_{l}}v_{\alpha }+\underset{\alpha ,p}{\sum }%
a_{\alpha }\frac{\partial ^{3}z^{\alpha }}{\partial z_{p}\partial
z_{k}\partial z_{l}}A_{p}^{j}+\underset{\alpha ,p}{\sum }a_{\alpha }\frac{%
\partial ^{3}z^{\alpha }}{\partial z_{j}\partial z_{p}\partial z_{l}}%
A_{p}^{k}+\underset{\alpha ,p}{\sum }a_{\alpha }\frac{\partial ^{3}z^{\alpha
}}{\partial z_{j}\partial z_{k}\partial z_{p}}A_{p}^{l}=0\text{ \ }(8_{jkl})
\end{equation*}

The equations for the unknowns $v_{\beta }^{\alpha }$ are obtained by
identifying the coefficients of the monomials $z^{\rho }$ in the above
equations. We can do the following reductions: using equations $6_{j}$ $%
v_{\beta }^{\alpha }=0$ if $\left| \alpha \right| +\left| \beta \right| \geq
d+1,$ and using the equation of $\mathcal{X}_{0}$, we can assume that degree
in the $z_{1}$ variable is at most $d-1.$

The monomials $z^{\rho }$ in (5) are $z_{1}^{\rho _{1}}$ $z_{2}^{\rho _{2}}$
$z_{3}^{\rho _{3}}$ $z_{4}^{\rho _{4}}$ with $\sum \rho _{i}\leq d$ and $%
\rho _{1}\leq d-1.$

If all the components of $\rho $ are greater than 3, then we obtain the
following system

9. The coefficient of $z^{\rho }$ in (5) impose the condition
\begin{equation*}
\underset{\alpha +\beta =\rho }{\sum }v_{\beta }^{\alpha }=0
\end{equation*}

10$_{j}.$ The coefficient of the monomial $z^{\rho -\delta _{j}}$ in $%
(6_{j}) $ impose the condition
\begin{equation*}
\underset{\alpha +\beta =\rho }{\sum }\alpha _{j}v_{\beta }^{\alpha
}=l_{j}(a)
\end{equation*}
where $l_{j}$ is a linear expression in the $a$-variables.

11$_{jj}.$ For $j=1,...,4$ the coefficient of the monomial $z^{\rho -2\delta
_{j}}$ in $(7_{jj})$ impose the condition
\begin{equation*}
\underset{\alpha +\beta =\rho }{\sum }\alpha _{j}(\alpha _{j}-1)v_{\beta
}^{\alpha }=l_{jj}(a)
\end{equation*}

11$_{jk}.$ For $1\leq j<k\leq 4$ the coefficient of the monomial $z^{\rho
-\delta _{j}-\delta _{k}}$ in $(7_{jk})$ impose the condition
\begin{equation*}
\underset{\alpha +\beta =\rho }{\sum }\alpha _{j}(\alpha _{j}-1)v_{\beta
}^{\alpha }=l_{jk}(a)
\end{equation*}

12$_{jjj}.$ For $j=1,...,4$ the coefficient of the monomial $z^{\rho
-3\delta _{j}}$ in $(8_{jjj})$ impose the condition
\begin{equation*}
\underset{\alpha +\beta =\rho }{\sum }\alpha _{j}(\alpha _{j}-1)(\alpha
_{j}-2)v_{\beta }^{\alpha }=l_{jjj}(a)
\end{equation*}

12$_{jjk}.$ For $1\leq j<k\leq 4$ the coefficient of the monomial $z^{\rho
-2\delta _{j}-\delta _{k}}$ in $(8_{jjk})$ impose the condition
\begin{equation*}
\underset{\alpha +\beta =\rho }{\sum }\alpha _{j}(\alpha _{j}-1)\alpha
_{k}v_{\beta }^{\alpha }=l_{jjk}(a)
\end{equation*}

12$_{jkl}.$ For $1\leq j<k<l\leq 4$ the coefficient of the monomial $z^{\rho
-\delta _{j}-\delta _{k}-\delta _{l}}$ in $(8_{jjk})$ impose the condition
\begin{equation*}
\underset{\alpha +\beta =\rho }{\sum }\alpha _{j}\alpha _{k}\alpha
_{l}v_{\beta }^{\alpha }=l_{jkl}(a)
\end{equation*}

As in the case of dimension 2 (cf. \cite{Paun05}) we obtain that the
determinant of the matrix associated to the system is not zero. Indeed, for
each $\rho $ the matrix whose column $C_{\beta }$ consists of the partial
derivatives of order at most 3 of the monomial $z^{\rho -\beta }$ has the
same determinant, at the point $z_{0}=(1,1,1,1),$ as our system. Therefore
if the determinant is zero, we would have a non-identically zero polynomial
\begin{equation*}
Q(z)=\underset{\beta }{\sum }a_{\beta }z^{\rho -\beta }
\end{equation*}
such that all its partial derivatives of order less or equal to 3 vanish at $%
z_{0}.$ Thus the same is true for
\begin{equation*}
P(z)=z^{\rho }Q(\frac{1}{z_{1}},...,\frac{1}{z_{4}})=\underset{\beta }{\sum }%
a_{\beta }z^{\beta }.
\end{equation*}
But this implies $P\equiv 0.$

Finally, we conclude by Cramer's rule. The systems we have to solve are
never over determined as well. The lemma is proved.
\end{proof}

\begin{proposition}
\label{p1}Let $\Sigma _{0}:=\{(z,a,\xi ^{(1)},\xi ^{(2)},\xi ^{(3)})\in
J_{3}^{v}(\mathcal{X})$ $/$ $\xi ^{(1)}\wedge \xi ^{(2)}\wedge \xi
^{(3)}=0\}.$ Then the vector space $T_{J_{3}^{v}(\mathcal{X)}}\otimes
\mathcal{O}_{\mathbb{P}^{4}}(12)\otimes \mathcal{O}_{\mathbb{P}%
^{N_{d}}}(\ast )$ is generated by its global sections on $J_{3}^{v}(\mathcal{%
X})\backslash \Sigma ,$ where $\Sigma $ is the closure of $\Sigma _{0}.$
\end{proposition}

\begin{proof}
From the preceding lemmas, we are reduced to consider $V=\underset{\left|
\alpha \right| \leq 3}{\sum }v_{\alpha }\frac{\partial }{\partial a_{\alpha }%
}.$ The conditions for $V$ to be tangent to $J_{3}^{v}(\mathcal{X}_{0})$ are
\begin{equation*}
\underset{\left| \alpha \right| \leq 3}{\sum }v_{\alpha }z^{\alpha }=0
\end{equation*}

\begin{equation*}
\underset{j=1}{\overset{4}{\sum }}\underset{\left| \alpha \right| \leq d,%
\text{ }\alpha _{1}<d}{\sum }v_{\alpha }\frac{\partial z^{\alpha }}{\partial
z_{j}}\xi _{j}^{(1)}=0
\end{equation*}

\begin{equation*}
\underset{\left| \alpha \right| \leq 3}{\sum }(\underset{j=1}{\overset{4}{%
\sum }}\frac{\partial z^{\alpha }}{\partial z_{j}}\xi _{j}^{(2)}+\underset{%
j,k=1}{\overset{4}{\sum }}\frac{\partial ^{2}z^{\alpha }}{\partial
z_{j}\partial z_{k}}\xi _{j}^{(1)}\xi _{k}^{(1)})v_{\alpha }=0
\end{equation*}

\begin{equation*}
\underset{\left| \alpha \right| \leq 3}{\sum }(\underset{j=1}{\overset{4}{%
\sum }}\frac{\partial z^{\alpha }}{\partial z_{j}}\xi _{j}^{(3)}+3\underset{%
j,k=1}{\overset{4}{\sum }}\frac{\partial ^{2}z^{\alpha }}{\partial
z_{j}\partial z_{k}}\xi _{j}^{(2)}\xi _{k}^{(1)}+\underset{j,k,l=1}{\overset{%
4}{\sum }}\frac{\partial ^{3}z^{\alpha }}{\partial z_{j}\partial
z_{k}\partial z_{l}}\xi _{j}^{(1)}\xi _{k}^{(1)}\xi _{l}^{(1)})v_{\alpha }=0
\end{equation*}

We denote by $W_{jkl}$ the wronskian operator corresponding to the variables
$z_{j},z_{k},z_{l}.$ We can suppose $W_{123}:=\det (\xi _{j}^{(i)})_{1\leq
i,j\leq 3}\neq 0.$ Then we can solve the previous system with $%
v_{0000},v_{1000},v_{0100},v_{0010}$ as unknowns. By the Cramer rule, each
of the previous quantity is a linear combination of the $v_{\alpha },$ $%
\left| \alpha \right| \leq 3,$ $\alpha \neq (0000),$ $(1000),$ $(0100),$ $%
(0010)$ with coefficients rational functions in $z,\xi ^{(1)},\xi ^{(2)},\xi
^{(3)}.$ The denominator is $W_{123}$ and the numerator is a polynomial
whose monomials verify either:

i) degree in $z$ at most 3 and degree in each $\xi ^{(i)}$ at most 1.

ii) degree in $z$ at most 2 and degree in $\xi ^{(1)}$ at most 3, degree in $%
\xi ^{(2)}$ at most 0, degree in $\xi ^{(3)}$ at most 1.

iii) degree in $z$ at most 2 and degree in $\xi ^{(1)}$ at most 2, degree in
$\xi ^{(2)}$ at most 2, degree in $\xi ^{(3)}$ at most 0.

iv) degree in $z$ at most 1 and degree in $\xi ^{(1)}$ at most 4, degree in $%
\xi ^{(2)}$ at most 1, degree in $\xi ^{(3)}$ at most 0.

$\xi ^{(1)}$ has a pole of order 2, $\xi ^{(2)}$ has a pole of order 3 and $%
\xi ^{(3)}$ has a pole of order 4, therefore the previous vector field has
order at most 12.
\end{proof}

\begin{remark}
\label{r1}If the third derivative of $f:(\mathbb{C},0)\rightarrow \mathcal{X}
$ lies inside $\Sigma _{0}$ then the image of $f$ is contained in a
hyperplane section of $\mathcal{X}.$
\end{remark}

\section{Jet differentials}

In this section we recall the basic facts about jet differentials following
J.-P. Demailly \cite{De95}.

Let $X$ be a complex manifold. We start with the directed manifold $%
(X,T_{X}) $. We define $X_{1}:=\mathbb{P(}T_{X}),$ and $V_{1}\subset
T_{X_{1}}:$%
\begin{equation*}
V_{1,(x,[v])}:=\{\xi \in T_{X_{1},(x,[v])}\text{ };\text{ }\pi _{\ast }\xi
\in \mathbb{C}v\}
\end{equation*}
where $\pi :X_{1}\rightarrow X$ is the natural projection. If $f:(\mathbb{C}%
,0)\rightarrow (X,x)$ is a germ of holomorphic curve then it can be lifted
to $X_{1}$ as $f_{[1]}.$

By induction, we obtain a tower of varieties $(X_{k},V_{k}).$ $\pi
_{k}:X_{k}\rightarrow X$ is the natural projection. We have a tautological
line bundle $\mathcal{O}_{X_{k}}(1)$ and we denote $u_{k}:=c_{1}(\mathcal{O}%
_{X_{k}}(1)).$

Let's consider the direct image $\pi _{k\ast }(\mathcal{O}_{X_{k}}(m)).$
It's a vector bundle over $X$ which can be described with local coordinates.
Let $z=(z_{1},...,z_{n})$ be local coordinates centered in $x\in X.$ A local
section of $\pi _{k\ast }(\mathcal{O}_{X_{k}}(m))$ is a polynomial
\begin{equation*}
P=\underset{\left| \alpha _{1}\right| +2\left| \alpha _{2}\right|
+...+k\left| \alpha _{k}\right| =m}{\sum }R_{\alpha }(z)dz^{\alpha
_{1}}...d^{k}z^{\alpha _{k}}
\end{equation*}
which acts naturally on the fibers of the bundle $J_{k}X\rightarrow X$ of $k$%
-jets of germs of curves in $X$, i.e the set of equivalence classes of
holomorphic maps $f:$ $(\mathbb{C},0)\rightarrow (X,x)$ with the equivalence
relation which identifies two such maps if their derivatives agree up to
order $k,$ and which is invariant under reparametrization i.e
\begin{equation*}
P((f\circ \phi )^{\prime },...,(f\circ \phi )^{(k)})_{t}=\phi ^{\prime
}(t)^{m}P(f^{\prime },...,f^{(k)})_{\phi (t)}
\end{equation*}
for every $\phi $ $\in $ $\mathbb{G}_{k},$ the group of $k$-jets of
biholomorphisms of $(\mathbb{C},0).$ The vector bundle $\pi _{k\ast }(%
\mathcal{O}_{X_{k}}(m))$ is denoted $E_{k,m}T_{X}^{\ast }.$ This bundle of
invariant jet differentials is a subbundle of the bundle of jet
differentials, of order $k$ and degree $m$, $E_{k,m}^{GG}T_{X}^{\ast
}\rightarrow X$ whose fibres are complex-valued polynomials $Q(f^{\prime
},f^{\prime \prime },...,f^{(k)})$ on the fibers of $J_{k}X,$ of weight $m$
under the action of $\mathbb{C}^{\ast }$:
\begin{equation*}
Q(\lambda f^{\prime },\lambda ^{2}f^{\prime \prime },...,\lambda
^{k}f^{(k)})=\lambda ^{m}Q(f^{\prime },f^{\prime \prime },...,f^{(k)})
\end{equation*}
for any $\lambda \in \mathbb{C}^{\ast }$ and $(f^{\prime },f^{\prime \prime
},...,f^{(k)})\in J_{k}X.$

It turns out that we have an embedding $J_{k}^{reg}X/\mathbb{G}%
_{k}\hookrightarrow X_{k},$ where $J_{k}^{reg}X$ denotes the space of
non-constants jets.

For $k=1,$ $E_{1,m}T_{X}^{\ast }=S^{m}T_{X}^{\ast }.$

If $X$ is a surface we have the following description of $E_{2,m}T_{X}^{\ast
}.$ Let $W$ be the wronskian, $W=dz_{1}d^{2}z_{2}-dz_{2}d^{2}z_{1},$ then
every invariant differential operator of order 2 and degree $m$ can be
written
\begin{equation*}
P=\underset{\left| \alpha \right| +3k=m}{\sum }R_{\alpha ,k}(z)dz^{\alpha
}W^{k}.
\end{equation*}
The following theorem makes clear the use of jet differentials in the study
of hyperbolicity:

\bigskip

\textbf{Theorem (\cite{GG80}, \cite{De95}). }\textit{Assume that there exist
integers }$k,m>0$\textit{\ and an ample line bundle }$L$\textit{\ on X such
that}
\begin{equation*}
H^{0}(X_{k},\mathcal{O}_{X_{k}}(m)\otimes \pi _{k}^{\ast }L^{-1})\simeq
H^{0}(X,E_{k,m}T_{X}^{\ast }\otimes L^{-1})
\end{equation*}
\textit{has non zero sections }$\sigma _{1},...,\sigma _{N}.$\textit{\ Let }$%
Z\subset X_{k}$\textit{\ be the base locus of these sections. Then every
entire curve }$f:\mathbb{C}\rightarrow X$\textit{\ tangent to }$V$\textit{\
is such that }$f_{[k]}(\mathbb{C})\subset Z.$\textit{\ In other words, for
every global }$\mathbb{G}_{k}-$\textit{invariant polynomial differential
operator P with values in }$L^{-1},$\textit{\ every entire curve }$f:\mathbb{%
C}\rightarrow X$\textit{\ tangent to V must satisfy the algebraic
differential equation }$P(f)=0.$

\bigskip

\begin{remark}
In fact, this theorem is true for global sections of $E_{k,m}^{GG}T_{X}^{%
\ast }$ vanishing on an ample divisor.
\end{remark}

A complex compact manifold is hyperbolic if there is no non constant entire
curve $f:\mathbb{C}\rightarrow X$. Thus, the problem reduces to produce
enough independant algebraic differential equations.

If $X\subset \mathbb{P}^{4}$ is a smooth hypersurface$,$ we have established
the next result:

\bigskip

\textbf{Theorem \cite{Rou2}. }\textit{Let X be a smooth hypersurface of }$%
\mathbb{P}^{4}$\textit{\ such that }$d=\deg (X)\geq 97,$\textit{\ and A an
ample line bundle, then }$E_{3,m}T_{X}^{\ast }\otimes A^{-1}$\textit{\ has
global sections for m large enough and every entire curve }$f:\mathbb{C}%
\rightarrow X$\textit{\ must satisfy the corresponding algebraic
differential equation.}

\bigskip

The proof relies on the filtration of $E_{3,m}T_{X}^{\ast }$ \cite{Rou05}$:$%
\begin{equation*}
Gr^{\bullet }E_{3,m}T_{X}^{\ast }=\underset{0\leq \gamma \leq \frac{m}{5}}{%
\oplus }(\underset{\{\lambda _{1}+2\lambda _{2}+3\lambda _{3}=m-\gamma ;%
\text{ }\lambda _{i}-\lambda _{j}\geq \gamma ,\text{ }i<j\}}{\oplus }\Gamma
^{(\lambda _{1},\lambda _{2},\lambda _{3})}T_{X}^{\ast })
\end{equation*}
where $\Gamma $ is the Schur functor. This filtration provides a
Riemann-Roch computation of the Euler characteristic \cite{Rou05}:
\begin{equation*}
\chi (X,E_{3,m}T_{X}^{\ast })=\frac{m^{9}}{81648\times 10^{6}}%
d(389d^{3}-20739d^{2}+185559d-358873)+O(m^{8}).
\end{equation*}
In dimension 3 there is no Bogomolov vanishing theorem (cf. \cite{Bo79}) as
it is used in dimension 2 to control the cohomology group $H^{2}$, therefore
we need the following proposition:

\bigskip

\textbf{Proposition \cite{Rou2}. }\textit{Let }$\lambda =(\lambda
_{1},\lambda _{2},\lambda _{3})$\textit{\ be a partition such that }$\lambda
_{1}>\lambda _{2}>\lambda _{3}$\textit{\ and }$\left| \lambda \right| =\sum
\lambda _{i}>4(d-5)+18.$\textit{\ Then : }
\begin{equation*}
h^{2}(X,\Gamma ^{\lambda }T_{X}^{\ast })\leq g(\lambda )d(d+13)+q(\lambda )
\end{equation*}
\textit{where }$g(\lambda )=\frac{3\left| \lambda \right| ^{3}}{2}\underset{%
\lambda _{i}>\lambda _{j}}{\prod }(\lambda _{i}-\lambda _{j})$\textit{\ and }%
$q$\textit{\ is a polynomial in }$\lambda $\textit{\ with homogeneous
components of degrees at most 5.}

\bigskip

This proposition provides the estimate
\begin{equation*}
h^{2}(X,Gr^{\bullet }E_{3,m}T_{X}^{\ast })\leq Cd(d+13)m^{9}+O(m^{8})
\end{equation*}
where C is a constant.

\section{Proof of the theorem}

Let us consider an entire curve $f:\mathbb{C}\rightarrow X$ in a generic
hypersurface of $\mathbb{P}^{4}.$ By Riemann-Roch and the proposition of the
previous section we obtain the following lemma:

\begin{lemma}
\textit{Let X be a smooth hypersurface of }$\mathbb{P}^{4}$ of degree $d$, $%
0<\delta <\frac{1}{18}$ then $h^{0}(X,E_{3,m}T_{X}^{\ast }\otimes
K_{X}^{-\delta m})\geq \alpha (d,\delta )m^{9}+O(m^{8}),$ with
\begin{eqnarray*}
\alpha (d,\delta ) &=&-\frac{1}{408240000000}d(-1945d^{3}-784080\delta
^{2}d^{3}+105030d^{3}\delta + \\
&&1058400\delta ^{3}d^{3}+103695d^{2}+7075491d+105837083+322256880\delta
^{2}d- \\
&&1260083250\delta -435002400\delta ^{3}d-6819271200\delta
^{3}+5051827440\delta ^{2}- \\
&&2255850d^{2}\delta -15876000\delta ^{3}d^{2}-81814050\delta
d+11761200\delta ^{2}d^{2}).
\end{eqnarray*}
\end{lemma}

\begin{proof}
$E_{3,m}T_{X}^{\ast }\otimes K_{X}^{-\delta m}$ admits a filtration with
graded pieces
\begin{equation*}
\Gamma ^{(\lambda _{1},\lambda _{2},\lambda _{3})}T_{X}^{\ast }\otimes
K_{X}^{-\delta m}=\Gamma ^{(\lambda _{1}-\delta m,\lambda _{2}-\delta
m,\lambda _{3}-\delta m)}T_{X}^{\ast }
\end{equation*}
for $\lambda _{1}+2\lambda _{2}+3\lambda _{3}=m-\gamma ;$ $\lambda
_{i}-\lambda _{j}\geq \gamma ,$ $i<j,$ $0\leq \gamma \leq \frac{m}{5}.$ We
compute by Riemann-Roch
\begin{equation*}
\chi (X,E_{3,m}T_{X}^{\ast }\otimes K_{X}^{-\delta m})=\chi (X,Gr^{\bullet
}E_{3,m}T_{X}^{\ast }\otimes K_{X}^{-\delta m}).
\end{equation*}
We use the proposition of the previous section to control
\begin{equation*}
h^{2}(X,E_{3,m}T_{X}^{\ast }\otimes K_{X}^{-\delta m}):
\end{equation*}
\begin{eqnarray*}
h^{2}(X,\Gamma ^{(\lambda _{1}-\delta m,\lambda _{2}-\delta m,\lambda
_{3}-\delta m)}T_{X}^{\ast }) &\leq &g(\lambda _{1}-\delta m,\lambda
_{2}-\delta m,\lambda _{3}-\delta m)d(d+13)+ \\
&&q(\lambda _{1}-\delta m,\lambda _{2}-\delta m,\lambda _{3}-\delta m)
\end{eqnarray*}
under the hypothesis $\sum \lambda _{i}-3\delta m>4(d-5)+18.$ The conditions
verified by $\lambda $ imply $\sum \lambda _{i}\geq \frac{m}{6}$ therefore
the hypothesis will be verified if
\begin{equation*}
m(\frac{1}{6}-3\delta )>4(d-5)+18.
\end{equation*}
We conclude with the computation
\begin{equation*}
\chi (X,E_{3,m}T_{X}^{\ast }\otimes K_{X}^{-\delta m})-h^{2}(X,Gr^{\bullet
}E_{3,m}T_{X}^{\ast }\otimes K_{X}^{-\delta m})\leq
h^{0}(X,E_{3,m}T_{X}^{\ast }\otimes K_{X}^{-\delta m}).
\end{equation*}
\end{proof}

\begin{remark}
If we denote $\mathcal{X}_{3}^{v}$ the quotient of $J_{3}^{v,reg}(\mathcal{X)%
}$ by the reparametrization group $\mathbb{G}_{3}$, one can easily verify
that each vector field given at section 2 defines a section of the tangent
bundle of the manifold $\mathcal{X}_{3}^{v}.$
\end{remark}

We have a section
\begin{equation*}
\sigma \in H^{0}(X,E_{3,m}T_{X}^{\ast }\otimes K_{X}^{-\delta m})\simeq
H^{0}(X_{3},\mathcal{O}_{X_{3}}(m)\otimes \pi _{3}^{\ast }K_{X}^{-\delta m}).
\end{equation*}
with zero set $Z$ and vanishing order $\delta m(d-5).$ Consider the family
\begin{equation*}
\mathcal{X}\subset \mathbb{P}^{4}\times \mathbb{P}^{N_{d}}
\end{equation*}
of hypersurfaces of degree $d$ in $\mathbb{P}^{4}.$ General semicontinuity
arguments concerning the cohomology groups show the existence of a Zariski
open set $U_{d}\subset \mathbb{P}^{N_{d}}$ such that for any $a\in U_{d},$
there exists a divisor
\begin{equation*}
Z_{a}=(P_{a}=0)\subset (\mathcal{X}_{a})_{3}
\end{equation*}
where
\begin{equation*}
P_{a}\in H^{0}((\mathcal{X}_{a})_{3},\mathcal{O}_{(\mathcal{X}%
_{a})_{3}}(m)\otimes \pi _{3}^{\ast }K_{(\mathcal{X}_{a})}^{-\delta m})
\end{equation*}
such that the family $(P_{a})_{a\in U_{d}}$ varies holomorphically. We
consider $P$ as a holomorphic function on $J_{3}(\mathcal{X}_{a}).$ The
vanishing order of this function as a function of $%
dz_{i},d^{2}z_{i},d^{3}z_{i}$ $(1\leq i\leq 3)$ is no more than $m$ at a
generic point of $\mathcal{X}_{a}.$ We have $f_{[3]}(\mathbb{C})\subset
Z_{a}.$

Then we invoke the proposition \ref{p1} which gives the global generation of
\begin{equation*}
T_{J_{3}^{v}(\mathcal{X)}}\otimes \mathcal{O}_{\mathbb{P}^{4}}(12)\otimes
\mathcal{O}_{\mathbb{P}^{N_{d}}}(\ast )
\end{equation*}
on $J_{3}^{v}(\mathcal{X})\backslash \Sigma .$

If $f_{[3]}(\mathbb{C)}$ lies in $\Sigma $, $f$ is algebraically degenerated
as we saw in remark \ref{r1}. So we can suppose it is not the case.

At any point of $f_{[3]}(\mathbb{C)}\backslash \Sigma $ where the vanishing
of $P$ as a function of $dz_{i},d^{2}z_{i},d^{3}z_{i}$ $(1\leq i\leq 3)$ is
no more than $m,$ we can find global meromorphic vector fields $%
v_{1},...,v_{p}$ $(p\leq m)$ and differentiate $P$ with these vector fields
such that $v_{1}...v_{p}P$ is not zero at this point$.$ From the above
remark, we see that $v_{1}...v_{p}P$ corresponds to an invariant
differential operator and its restriction to $(\mathcal{X}_{a})_{3}$ can be
seen as a section of the bundle
\begin{equation*}
\mathcal{O}_{(\mathcal{X}_{a})_{3}}(m)\otimes \mathcal{O}_{\mathbb{P}%
^{4}}(12p-\delta m(d-5)).
\end{equation*}
Assume that the vanishing order of $P$ is larger than the sum of the pole
order of the $v_{i}$ in the fiber direction of $\pi :\mathcal{X}\rightarrow
\mathbb{P}^{N_{d}}.$ Then the restriction of $v_{1}...v_{p}P$ to $\mathcal{X}%
_{a}$ defines a jet differential which vanishes on an ample divisor.
Therefore $f_{[3]}(\mathbb{C)}$ should be in its zero set. Thus $f_{[3]}(%
\mathbb{C)}$ should be in the zero section of $J_{3}(\mathcal{X}_{a})$ over
a generic point of $\mathcal{X}_{a}.$

To finish the proof, we just have to see when the vanishing order of $P$ is
larger than the sum of the pole order of the $v_{i}.$ This will be verified
if
\begin{equation*}
\delta (d-5)>12.
\end{equation*}
So we want $\delta >\frac{12}{(d-5)}$ and $\alpha (d,\delta )>0.$ This is
the case for $d\geq 593.$

\end{document}